\documentclass[12pt,draft,leqno]{article}
\usepackage{amssymb, eucal, latexsym}
\newtheorem{theorem}{Theorem}[section]
\newtheorem{lm}[theorem]{Lemma}

\newtheorem{cor}[theorem]{Corollary}
\newtheorem{pro}[theorem]{Proposition}
\newtheorem{defi}[theorem]{Definition}

\newtheorem{nota}[theorem]{Notation}
\newtheorem{notas}[theorem]{Notations}
\newtheorem{rem}[theorem]{Remark}

\newtheorem{fact}[theorem]{Fact}

\def\p{\varphi}
\def\a{\alpha}
\def\b{\beta}
\def\d{\delta}

\def\l{\lambda}

\def\s{\sigma}

\def\OM{\Omega}

\def\pl{\varphi_\Lambda}
\def\fs{\hat{f}}
\def\ps{\hat{\varphi}}

\def\lra{\longrightarrow}

\def\sbe{\subseteq}
\def\spe{\supseteq}
\def\stm{\setminus}
\def\ems{\emptyset}
\def\nes{\neq\emptyset}

\def\ex{\exists}
\def\fa{\forall}

\def\we{\wedge}

\def\ap{^\prime}
\def\inv{^{-1}}
\def\st{\ |\ }

\def\llx{\ll_{\rho}}

\def\nin{\not\in}

\def\BB{{\cal B}}

\def\KK{{\cal K}}
\def\LL{{\cal L}}

\def\OO{{\cal O}}
\def\PP{{\cal P}}

\def\ISKLC{{\bf InSkeLC}}
\def\SSKLC{{\bf SuSkeLC}}
\def\ISKAL{{\bf DSuSkeLC}}
\def\ISKLC{{\bf InSkeLC}}
\def\SSKAL{{\bf DInSkeLC}}
\def\IOPLC{{\bf InOpPerLC}}
\def\SOPLC{{\bf SuOpPerLC}}
\def\IOPAL{{\bf DSuOpPerLC}}
\def\SOPAL{{\bf DInOpPerLC}}

\def\IOHC{{\bf InOpC}}
\def\IOLC{{\bf InOpLC}}
\def\IOAC{{\bf DSuOpC}}
\def\IOAL{{\bf DSuOpLC}}
\def\SOHC{{\bf SuOpC}}
\def\SOLC{{\bf SuOpLC}}
\def\SOAC{{\bf DInOpC}}
\def\SOAL{{\bf DInOpLC}}

\def\SAC{{\bf DSkeC}}
\def\SAL{{\bf DSkePerLC}}
\def\SKAL{{\bf DSkeLC}}
\def\SKLC{{\bf SkeLC}}
\def\ISHC{{\bf InSkeC}}
\def\ISLC{{\bf InSkePerLC}}
\def\ISAC{{\bf DSuSkeC}}
\def\ISAL{{\bf DSuSkePerLC}}
\def\SSHC{{\bf SuSkeC}}
\def\SSLC{{\bf SuSkePerLC}}
\def\SSAC{{\bf DInSkeC}}
\def\SSAL{{\bf DInSkePerLC}}

\def\2{\mbox{{\bf 2}}}
\def\3{\mbox{{\bf 3}}}

\def\int{\mbox{{\rm int}}}

\def\cl{\mbox{{\rm cl}}}

\def\Ult{\mbox{{\rm Ult}}}

\def\doc{\hspace{-1cm}{\em Proof.}~~}
\def\sq{\hspace*{\fill} \hbox{\vrule\vbox{\hrule\phantom{o}\hrule}\vrule}}
\def\sqs{\sq \vspace{2mm}}

\def\cp{compactification }

\def\BBBB{{\rm I}\!{\rm B}}
\def\DDDD{{\rm I}\!{\rm D}}

\title{{\LARGE\bf
Some Generalizations of Fedorchuk Duality Theorem -- II}\\
\vspace{0.35cm}
{\large\bf Georgi Dimov}\thanks{This paper was supported by the
project no. 101/2007 $``$Categorical Topology" of the Sofia
University $``$St. Kl. Ohridski".}\\
\vspace{0.25cm}
 {\footnotesize Dept. of Math. and
Informatics, Sofia University,  Blvd. J. Bourchier 5, 1164 Sofia,
Bulgaria}}

\author{}

\date{}

\begin{document}
\maketitle
\begin{abstract}
{\footnotesize
\noindent This paper is a continuation of the paper \cite{Di2}. In
\cite{Di2} it was shown that there exists a duality $\Psi^a$
between the category $\SKAL$ (introduced there) and the category
$\SKLC$ of locally compact Hausdorff spaces and continuous
skeletal maps. We describe here the subcategories of the category
$\SKAL$ which are dually equivalent to the following eight
categories: all of them have as objects the locally compact
Hausdorff spaces and their morphisms are, respectively, the
injective (respectively, surjective) continuous skeletal maps, the
injective (surjective) open maps, the injective (surjective)
skeletal perfect maps,  the injective (surjective) open perfect
maps.  The particular cases of these theorems for the full
subcategories of the last four categories having as objects all
compact Hausdorff spaces are formulated and proved. The
$\SKAL$-morphisms which are LCA-embeddings and the dense
homeomorphic embeddings are characterized through their dual morphisms.
For any locally compact space $X$, a description of the frame of
all open subsets of $X$ in terms of the dual object of $X$ is
obtained. It is shown how one can build the dual object of an open
subset (respectively, of a regular closed subset) of a locally
compact Hausdorff space $X$ directly from the dual object of $X$.
Applying these results, a new description of the ordered set of
all, up to equivalence, locally compact Hausdorff extensions of a
locally compact Hausdorff space is obtained. Moreover,
generalizing de Vries Compactification Theorem (\cite{dV2}), we
strengthen the Local Compactification Theorem of Leader
(\cite{LE2}). Some other applications are found.}
\end{abstract}

{\footnotesize {\em  MSC:} primary 54D45, 18A40; secondary 54C10,
54D35, 54E05.

}

\footnotetext[1]{{\footnotesize {\em E-mail address:}
gdimov@fmi.uni-sofia.bg}}

\baselineskip = \normalbaselineskip

\section*{Introduction}

This paper is the second part of the paper \cite{Di2}. In
\cite{Di2},  the category $\SKAL$ of all complete local contact
algebras and all complete Boolean homomorphisms between them
satisfying two simple conditions (see \cite[Definition 2.10]{Di2})
was defined and it was shown that there exists a duality
$\Psi^a:\SKAL\lra\SKLC$, where $\SKLC$ is the category of all
locally compact Hausdorff spaces and all skeletal (in the sense of
Mioduszewski and Rudolf \cite{MR2}) continuous maps between them.
In the first section of the present paper, we find the
subcategories of the category $\SKAL$ which are dually equivalent
to the following eight subcategories of the category $\SKLC$: all
of them have as objects the locally compact Hausdorff spaces and
their morphisms are, respectively, the injective (respectively,
surjective) continuous skeletal maps (see Theorem \ref{ilcompnk}
(resp., Theorem \ref{slcompnk})), the injective (surjective) open
maps (see Theorem \ref{ilcompo} (resp., Theorem \ref{slcompo})),
the injective (surjective) skeletal perfect maps (see Theorem
\ref{ilcompn} (resp., Theorem \ref{slcompn})), the injective
(surjective) open perfect maps (see Theorem \ref{ilcomp} (resp.,
Theorem \ref{slcomp})). In the theorems mentioned above, the
particular cases for the full subcategories of the last four
categories having as objects all compact Hausdorff spaces are
formulated as well.

In the second section, we prove that $\p$ is a $\SKAL$-morphism
and an LCA-embedding iff $f=\Psi^a(\p)$ is a quasi-open semi-open
perfect surjection (see Theorem \ref{lcaemb}). The dense
homeomorphic embeddings are characterized  through their dual morphisms
as well (see Theorem \ref{f1}). For any locally compact Hausdorff
space $X$, a description of the frame of all open subsets of $X$
in terms of the dual to $X$ complete local contact algebra is
obtained (see Theorem \ref{opensetsfr}). It is shown how one can
build the dual object of an open subset (respectively, of a
regular closed subset) of a locally compact Hausdorff space $X$
directly from the CLCA dual to $X$ (see Theorem \ref{conopen}
(respectively, Theorem \ref{conregclo})).

The third section is devoted to some immediate applications of the
theorems obtained above. A new description of the ordered set of
all (up to equivalence) Hausdorff locally compact extensions of a
locally compact Hausdorff space is obtained in the language of
normal contact relations (see Theorem \ref{compset}). The normal
contact relations which correspond to the Alexandroff (one-point)
compactification and to the Stone-\v{C}ech compactification of a
locally compact Hausdorff space are found (see Theorem
\ref{alphabeta}). The Wallman-type compactifications of discrete
spaces are described as well (see Theorem \ref{Wallmantype}).
Generalizing de Vries Compactification Theorem (\cite{dV2}) (see
Corollary \ref{lead7}), a theorem is  obtained (see Theorem
\ref{lead5}) which  strengthen the Local Compactification Theorem
of Leader (\cite{LE2}). Some other applications are found.

In this paper we use the notations introduced in the first part of
it (see \cite{Di2}), as well as the notions and the results from
\cite{Di2}. If $A$ is a Boolean algebra then we denote by
$Atoms(A)$ the set of all atoms of $A$. If $(A,\le)$ is a poset
and $a\in A$ then  $\downarrow a$ is the set $\{b\in A\st b\le
a\}$. If $f:X\lra Y$ is a function and $M\sbe X$ then
$f\upharpoonright M$ is the restriction of $f$ having $M$ as a
domain and $f(M)$ as a codomain. Finally, we will denote by
$\DDDD$  the set of all dyadic numbers of the interval $(0,1)$.

\section{Surjective and Injective Mappings}

\begin{notas}\label{newcategnk}
\rm We  denote by:

\noindent $\bullet$ $\ISKLC$ the category of all locally compact
Hausdorff spaces and all continuous skeletal injective  maps
between them;

\noindent $\bullet$ $\SSKLC$ the category of all locally compact
Hausdorff spaces and all continuous skeletal surjective maps
between them.

 In the sequel we will use similar notations without more explanations, i.e.
  if {\bf K} is a  category introduced in \cite{Di2},
 then  by {\bf InK} (resp., {\bf SuK}) we will denote the
category  having the same objects as the category {\bf K} and
whose morphisms are only the injective (resp., surjective)
morphisms of {\bf K}.
\end{notas}

\begin{defi}\label{surk}
\rm
Let $\ISKAL$ be the category whose objects are all complete local
contact algebras (see \cite[1.14]{Di2}) and whose morphisms are
all $\SKAL$-morphisms (see \cite[2.10]{Di2})
$\p:(A,\rho,\BBBB)\lra (B,\eta,\BBBB\ap)$
 which satisfy the following condition:

\smallskip

\noindent (IS) For every bounded ultrafilter $u$ in
$(A,\rho,\BBBB)$ there exists a bounded ultrafilter $v$ in
$(B,\eta,\BBBB\ap)$ such that $\pl(v)\rho u$ (see \cite[1.18]{Di2}
and \cite[1.21]{Di2} for the notations).

\smallskip

Using \cite[Corollary 1.5, (EL1), (36), 1.21]{Di2},
%\cite[(\ref{ultbasis})]{Di2}, \cite[\ref{ladj}]{Di2},
we obtain easily that $\ISKAL$ is indeed a category.
\end{defi}

\begin{theorem}\label{slcompnk}
The categories $\SSKLC$ and $\ISKAL$ are dually equivalent.
\end{theorem}

\doc Let $f:X\lra Y$ be a surjective continuous skeletal map
between two locally compact Hausdorff spaces and $\p=\Psi^t(f)$.
Then $\p:RC(Y)\lra RC(X)$ and $\pl(F)=\cl(f(F))$, for every $F\in
RC(X)$ (see the proof of \cite[Theorem 2.11]{Di2}). Let $u$ be a
bounded ultrafilter in $RC(Y)$. Then there exists $G_0\in
CR(Y)\cap u$. Hence there exists $y\in\bigcap\{G\st G\in u\}$.
Since $f$ is a surjection, there exists $x\in X$ such that
$f(x)=y$. Let $v$ be an ultrafilter in $RC(X)$ which contains
$\nu_x$ (see \cite[(3)]{Di2} for $\nu_x$). Then, obviously, $v$ is
a bounded ultrafilter in $(RC(X),\rho_X,CR(X))$. By
\cite[(51)]{Di2}, $v\sbe\s_x$ (see \cite[(3)]{Di2} for $\s_x$).
Hence $y\in\pl(F)$, for every $F\in v$. This means that
$\pl(v)\rho_Y u$. Therefore, $\p$ satisfies condition (IS). By
\cite[Theorem 2.11]{Di2}, $\p$ is a $\SKAL$-morphism. Hence, $\p$
is a $\ISKAL$-morphism.

Let $\p:(A,\rho,\BBBB)\lra (B,\eta,\BBBB\ap)$ be a
$\ISKAL$-morphism and $f=\Psi^a(\p)$. Let
$X=\Psi^a(A,\rho,\BBBB)$, $Y=\Psi^a (B,\eta,\BBBB\ap)$ and $\s\in
X$. Then $\s$ is a bounded cluster in $(A,C_\rho)$. Hence there
exists a bounded ultrafilter $u$ in $(A,\rho,\BBBB)$ such that
$\s=\s_u$. By (IS), there exists a bounded ultrafilter $v$ in
$(B,\eta,\BBBB\ap)$ such that $\pl(v)\rho u$. Thus $\pl(v)C_\rho
u$. Therefore, by \cite[1.5, (36), (37)]{Di2},
%\cite[(\ref{ultbasis})]{Di2}, \cite[(\ref{plu})]{Di2}
and \cite[(35)]{Di2}, $f(\s_v)=\s_{\pl(v)}=\s_u=\s$. So, $f$ is a
surjection. By  \cite[Theorem 2.11]{Di2}, $f$ is a continuous
skeletal map. Hence,  $f$ is a $\SSKLC$-morphism.

The rest follows from  \cite[Theorem 2.11]{Di2}. \sqs

\begin{cor}\label{inj}
Every $\ISKAL$-morphism is an injection.
\end{cor}

\doc Let $\p:(A,\rho,\BBBB)\lra
(B,\eta,\BBBB\ap)$ be a $\ISKAL$-morphism. Set $f=\Psi^a(\p)$,
$X=\Psi^a(B,\eta,\BBBB\ap)$ and $Y=\Psi^a(A,\rho,\BBBB)$. Then, by
Theorem \ref{slcompnk}, $f:X\lra Y$ is a surjective continuous
skeletal map.

Let $\ps=\Psi^t(f)$. Then $\ps:RC(Y)\lra RC(X)$. We will show that
$\ps$ is an injection. We have that, for every $F\in RC(Y)$,
$\ps(F)=\cl(\int(f\inv(F)))=\cl(f\inv(\int(F)))$. Let $F,G\in
RC(Y)$ and $\ps(F)=\ps(G)$. Suppose that there exists a point
$x\in F\stm G$. Then there exists an open neighborhood $Ox$ of $x$
such that $Ox\cap G=\ems$. Obviously, there exists a point $y\in
Ox\cap \int(F)$. Then $f\inv(y)\cap f\inv(G)=\ems$. But
$f\inv(y)\sbe f\inv(\int(F))\sbe\ps(F)=\ps(G)\sbe f\inv(G)$ and
thus $f\inv(y)=\ems$. Since $f$ is a surjection, we obtain a
contradiction. Hence $F\sbe G$. Analogously we prove that $G\sbe
F$. Therefore $F=G$. So, $\ps$ is an injection. Now, using Theorem
\ref{slcompnk}, we get that $\p$ is an injection.
\sqs

\begin{defi}\label{soamorphik}
\rm
Let $\SSKAL$ be the category whose objects are all complete local
contact algebras (see \cite[1.14]{Di2}) and whose morphisms are
all $\SKAL$-morphisms (see \cite[2.10]{Di2})
$\p:(A,\rho,\BBBB)\lra (B,\eta,\BBBB\ap)$ which satisfy the
following condition:

\smallskip

\noindent (LS) $\fa a,b\in \BBBB\ap$, $\pl(a)\rho\pl(b)$ implies
$a\eta b$ (see \cite[1.21]{Di2} for $\pl$).

\smallskip

It is easy to see that $\SSKAL$ is indeed a category.
\end{defi}

\begin{theorem}\label{ilcompnk}
The categories $\ISKLC$ and $\SSKAL$ are dually equivalent.
\end{theorem}

\doc Let $f:X\lra Y$ be an injective continuous skeletal map
between two locally compact Hausdorff spaces $X$ and $Y$. Set
$\p=\Psi^t(f)$ (see \cite[(31)]{Di2} for $\Psi^t(f)$). The
function $\pl:RC(X)\lra RC(Y)$ is defined by $\pl(F)=\cl(f(F))$,
for every $F\in RC(X)$ (see \cite[(32) and (33)]{Di2}). Hence, for
$F\in CR(X)$, $\pl(F)=f(F)$. Since $f$ is an injection, it becomes
obvious that $\p$ satisfies condition (LS) from \ref{soamorphik}.
By \cite[Theorem 2.11]{Di2}, $\Psi^t(f)$ is a $\SKAL$-morphism.
Thus we get that $\Psi^t(f)$ is a $\SSKAL$-morphism.

Let $\p:(A,\rho,\BBBB)\lra (B,\eta,\BBBB\ap)$ be a morphism of the
category $\SSKAL$ and let $f=\Psi^a(\p)$ (see \cite[(35)]{Di2} for
$\Psi^a(\p)$). We will show that $f$ is an injection. Let $\s$ and
$\s\ap$ be two bounded clusters in $(B,C_\eta)$ (see
\cite[1.16]{Di2} for $C_\eta$) and $\s\neq\s\ap$. Since
$Y=\Psi^a(B,\eta,\BBBB\ap)$ is a locally compact Hausdorff space,
we obtain, using \cite[(24)]{Di2}, that there exist
$b,b\,\ap\in\BBBB\ap$ such that $\s\in\int(\l_B^g(b))$,
$\s\ap\in\int(\l_B^g(b\,\ap))$ and
$\l_B^g(b)\cap\l_B^g(b\,\ap)=\ems$. Thus, by \cite[2.1]{Di2} and
\cite[(14)]{Di2}, $b\,(-\eta)b\,\ap$.  Now, (LS) implies that
$\pl(b)(-\rho)\pl(b\,\ap)$. Using condition (L2) from \cite{Di2},
we get that $\pl(b),\pl(b\,\ap)\in\BBBB$. Hence
$\pl(b)(-C_\rho)\pl(b\,\ap)$. Since, by \cite[(40)]{Di2},
$\pl(b)\in f(\s)$ and $\pl(b\,\ap)\in f(\s\ap)$, we obtain that
$f(\s)\neq f(\s\ap)$. So, $f$ is an injection. Finally, by
\cite[Theorem 2.11]{Di2}, $f$ is a continuous skeletal map. Thus
$f$ is an $\ISKLC$-morphism.

Now,  \cite[Theorem 2.11]{Di2} implies that the categories
$\ISKLC$ and $\SSKAL$ are dually equivalent.
\sqs

\begin{notas}\label{newcategn}
\rm We  denote by:

\noindent $\bullet$ $\ISAL$ the category of all CLCA's and all
injective complete Boolean homomorphisms between them satisfying
axioms (L1)-(L3) (see \cite[Definitions 2.10, 2.14]{Di2});

\noindent $\bullet$ $\ISAC$ the category of all CNCA's and all
injective complete Boolean homomorphisms between them satisfying
axiom (F1) (see \cite[Definition 2.12]{Di2}).
\end{notas}

\begin{theorem}\label{slcompn}
The categories $\SSLC$ and $\ISAL$ are dually equivalent; in
particular, the categories $\SSHC$ and $\ISAC$ are dually
equivalent.
\end{theorem}

\doc Let $f:X\lra Y$ be a surjective skeletal perfect map between
two locally compact Hausdorff spaces $X$ and $Y$. Then the proof
of Corollary \ref{inj} shows that $\Psi^t(f):RC(Y)\lra RC(X)$ is
an injection. By  \cite[Theorem 2.15]{Di2}, $\Psi^t(f)$ is an
$\SAL$-morphism. So, we get that $\Psi^t(f)$ is a
$\ISAL$-morphism.

Let $\p:(A,\rho,\BBBB)\lra (B,\eta,\BBBB\ap)$ be a
$\ISAL$-morphism. Then $\p$ is an injection. Hence $\pl(\p(a))=a$,
for every $a\in A$. We will show that $\p$ satisfies condition
(IS) from \ref{surk}. Let $u$ be a bounded ultrafilter in
$(A,\rho,\BBBB)$. Then there exists $a\in u\cap\BBBB$. Since $\p$
is an injection, $\p(u)$ is a filter base in $B$. Hence, there
exists $v\in\Ult(B)$ such that $\p(u)\sbe v$. By the condition
(L3) from \cite{Di2}, $\p(a)\in\BBBB\ap$. Therefore, $v$ is a
bounded ultrafilter in $(B,\eta,\BBBB\ap)$.  Moreover, $\pl(v)=u$.
Indeed, since $\p(u)\sbe v$ and $\pl\circ\p=id_A$, we get that
$u\sbe\pl(v)$. Further, using \cite[1.21]{Di2}, we obtain easily
that $\pl(v)$ is a filter base.
 Hence, $\pl(v)=u$. Therefore $\pl(v)\rho u$. Thus
$\p$ satisfies condition (IS). Then Theorem \ref{slcompnk} implies
that $f=\Psi^a(\p)$ is a surjection. By \cite[Theorem 2.15]{Di2},
$f$ is a skeletal perfect map. So, we get that $f$ is a
$\SSLC$-morphism.

Now, using  \cite[Theorem 2.15]{Di2}, we obtain that the
categories $\SSLC$ and $\ISAL$ are dually equivalent. In
particular, the categories $\SSHC$ and $\ISAC$ are dually
equivalent.
\sqs

\begin{defi}\label{soamorphi}
\rm
Let $\SSAL$ be the category whose objects are all complete local
contact algebras (see \cite[1.14]{Di2}) and whose morphisms are
all $\SAL$-morphisms (see \cite[2.14]{Di2}) which satisfy
condition (LS) (see \ref{soamorphik}).

\smallskip

Let $\SSAC$ be the category whose objects are all complete normal
contact algebras (see \cite[1.1]{Di2}) and whose morphisms are all
$\SAC$-morphisms (see \cite[2.12]{Di2}) $\p:(A,C)\lra (B,C\ap)$
which satisfy the following condition:

\smallskip

\noindent (FS) $\fa a,b\in B$, $\pl(a)C\pl(b)$ implies $aC\ap b$.
\end{defi}

\begin{theorem}\label{ilcompn}
The categories $\ISLC$ and $\SSAL$ are dually equivalent; in
particular, the categories $\ISHC$ and $\SSAC$ are dually
equivalent.
\end{theorem}

\doc It follows from  \cite[Theorem 2.15]{Di2} and Theorem
\ref{ilcompnk}.
\sqs

Note that since the morphisms of the category $\ISLC$ are closed
maps, in the definition of the category $\SSAL$ (see
\ref{soamorphi}) we can substitute condition (LS)  for the
following one:

\smallskip

\noindent (ELS) $\fa a,b\in B$, $\pl(a)\rho\pl(b)$ implies $a\eta
b$ (see \cite[1.21]{Di2} for $\pl$).

\begin{notas}\label{newcatego}
\rm We  denote by:

\noindent $\bullet$ $\IOAL$ the category of all CLCA's and all
complete Boolean homomorphisms between them satisfying axioms
(L1), (L2),  (IS) and (LO) (see \cite[2.10]{Di2}, \ref{surk} and
\cite[2.16]{Di2});

\noindent $\bullet$ $\SOAL$ the category of all CLCA's and all
surjective complete Boolean homomorphisms between them satisfying
axioms (L1), (L2) and (LO) (see \cite[2.10]{Di2} and
\cite[2.16]{Di2}).
\end{notas}

\begin{theorem}\label{ilcompo}
The categories $\IOLC$ and $\SOAL$ are dually equivalent.
\end{theorem}

\doc Let us show that every $\SOAL$-morphism is a $\SSKAL$-morphism.
Indeed, let $\p:(A,\rho,\BBBB)\lra (B,\eta,\BBBB\ap)$ be a
morphism of the category $\SOAL$. Then $\p$ is a surjection. Let
$a,b\in \BBBB\ap$ and $\pl(a)\rho\pl(b)$. Then, by condition (LO)
 from \cite{Di2}, $\p(\pl(a))\eta b$. Hence, by
\cite[1.21]{Di2}, $a\eta b$. Therefore, $\p$ satisfies condition
(LS). Hence $\p$ is a $\SSKAL$-morphism.

Let $f:X\lra Y$ be an $\IOLC$-morphism and $\p=\Psi^t(f)$. Then
$\p(G)=f\inv(G)$, for every $G\in RC(Y)$ (see  the proof of
\cite[Theorem 2.17]{Di2}). For every $F\in RC(X)$ we have, by
\cite[Corollary 2.5]{Di2} and \cite[Lemma 2.6]{Di2}, that
$\cl(f(F))\in RC(Y)$. Set $G=\cl(f(F))$. Then, by
\cite[1.4.C]{E2}, $f\inv(G)=\cl(f\inv(f(F)))$ (because $f$ is an
open map), and the injectivity of $f$ implies that $f\inv(G)=F$.
Hence $\p(G)=F$. Therefore, $\p$ is a surjection.

 Now, our
theorem follows from Theorem \ref{ilcompnk} and \cite[Theorem
2.17]{Di2}.
\sqs

\begin{theorem}\label{slcompo}
The categories $\SOLC$ and $\IOAL$ are dually equivalent.
\end{theorem}

\doc It follows from Theorem \ref{slcompnk} and
\cite[Theorem 2.17]{Di2}.
\sqs

\begin{notas}\label{newcateg}
\rm We will denote by:

\noindent $\bullet$ $\IOPAL$ the category of all CLCA's and all
injective complete Boolean homomorphisms between them satisfying
axioms (L1)-(L3) and (LO) (see \cite[Definitions 2.10, 2.14,
2.16]{Di2});

\noindent $\bullet$ $\SOPAL$ the category of all CLCA's and all
surjective complete Boolean homomorphisms between them satisfying
axioms (L1)-(L3) and (LO) (see \cite[Definitions 2.10, 2.14,
2.16]{Di2});

\noindent $\bullet$ $\IOAC$ the category of all CNCA's and all
injective complete Boolean homomorphisms between them satisfying
axioms (CO) and (F1)  (see \cite[Definitions 2.12, 2.18]{Di2});

\noindent $\bullet$ $\SOAC$ the category of all CNCA's and all
surjective complete Boolean homomorphisms between them satisfying
axioms (CO) and (F1) (see \cite[Definitions 2.12, 2.18]{Di2}).
\end{notas}

\begin{theorem}\label{ilcomp}
The categories $\IOPLC$ and $\SOPAL$ are dually equivalent; in
particular, the categories $\IOHC$ and $\SOAC$ are dually
equivalent.
\end{theorem}

\doc It follows from  Theorem
\ref{ilcompo} and  \cite[Theorem 2.21]{Di2}.
\sqs

\begin{theorem}\label{slcomp}
The categories $\SOPLC$ and $\IOPAL$ are dually equivalent; in
particular, the categories $\SOHC$ and $\IOAC$ are dually
equivalent.
\end{theorem}

\doc It follows from Theorem \ref{slcompn} and
\cite[Theorem 2.21]{Di2}.
\sqs

\section{Embeddings. Open sets. Regular closed sets}

 We will need a lemma from \cite{CNG}:

\begin{lm}\label{isombool}
Let $X$ be a dense subspace of a topological space $Y$. Then the
functions $r:RC(Y)\lra RC(X)$, $F\lra F\cap X$, and $e:RC(X)\lra
RC(Y)$, $G\lra \cl_Y(G)$, are Boolean isomorphisms between Boolean
algebras $RC(X)$ and $RC(Y)$, and $e\circ r=id_{RC(Y)}$, $r\circ
e=id_{RC(X)}$.
\end{lm}

\begin{theorem}\label{f1}
If a function $f$ between two locally compact Hausdorff spaces is
a dense homeomorphic embedding then $\p=\Psi^t(f)$ is a
$\SOAL$-morphism and a  Boolean isomorphism. Conversely, if a
function $\p$ is  a $\SOAL$-morphism and a Boolean isomorphism
then $f=\Psi^a(\p)$ is a dense homeomorphic embedding.
\end{theorem}

\doc Let $f$ be a dense homeomorphic embedding of $X$ in $Y$. Then
$f(X)$ is a locally compact dense subspace of $Y$ and hence it is
open in $Y$. Thus $f$ is an open injection, i.e.  $f$ is an
$\IOLC$-morphism. Therefore, by Theorem \ref{ilcompo}, $\p$ is a
$\SOAL$-morphism. Put $Z=f(X)$ and let $i:Z\lra Y$ be the
embedding of $Z$ in $Y$. Then $\psi=\Psi^t(i):RC(Y)\lra RC(Z)$ is
defined by the formula $\psi(F)=cl_Z(Z\cap\int_Y(F))=F\cap Z$, for
every $F\in RC(Y)$. Hence, by Lemma \ref{isombool}, $\psi$ is a
Boolean isomorphism. Since $f=i\circ (f\upharpoonright X)$, we
obtain that $\p$ is a Boolean isomorphism as well.

Conversely, if $\p:(A,\rho,\BBBB)\lra(B,\eta,\BBBB\ap)$ is a
$\SOAL$-morphism and a  Boolean isomorphism  then, by Theorem
\ref{ilcompo}, $f=\Psi^a(\p)$ is a homeomorphic embedding. Let
$X=\Psi^a(B,\eta,\BBBB\ap)$, $Y=\Psi^a(A,\rho,\BBBB)$ and
$\ps=\Psi^t(f)$. Then Theorem \ref{ilcompo} shows that $\ps$ is a
Boolean isomorphism. Thus, if $F\in RC(Y)$ and $F\nes$ then
$\ps(F)\nes$, i.e. $\cl(f\inv(\int(F))\nes$. This means that
$f(X)\cap\int(F)\nes$. Therefore, $f(X)$ is dense in $Y$.
\sqs

\begin{pro}\label{genemb}
Let $f\in\SKLC(X,Y)$ and $\p=\Psi^t(f)$. Then $f$ is a
homeomorphic embedding iff $\p=\p_1\circ\p_2$ where  $\p_1$ is a
$\SSAL$-morphism and $\p_2$ is a Boolean isomorphism and a
$\SOAL$-morphism.
\end{pro}

\doc Let $f:X\lra Y$ be a homeomorphic embedding. Put
$Z=cl_Y(f(X))$ and let $f_1:X\lra Z$ be the restriction of $f$,
$f_2:Z\lra Y$ be the embedding of $Z$ in $Y$. Then $f=f_2\circ
f_1$. Obviously, $f_1$ is an open map. Let $U$ be an open subset
of $Z$. Then $V=f\inv(U)$ is open in $X$ and $f(V)\sbe f_2(U)$.
Now, the fact that $f$ is a skeletal map and \cite[Lemma 2.4]{Di2}
imply that $f_2$ is a skeletal map. Hence, $f_i$, $i=1,2$, are
$\SKLC$-morphisms. Set $\p_i=\Psi^t(f_i)$, $i=1,2$. Then
$\p=\p_1\circ\p_2$ and Theorem \ref{f1} together with Theorem
\ref{ilcompn} show that $\p_i$, $i=1,2$, are as required.

Conversely, let $\p=\p_1\circ\p_2$, where   $\p_1$ is a
$\SSAL$-morphism and $\p_2$ is a $\SOAL$-morphism and a Boolean
isomorphism. Set $\fs=\Psi^a(\p)$ and $\fs_i=\Psi^a(\p_i)$,
$i=1,2$. Then $\fs=\fs_2\circ \fs_1$ and, by Theorems \ref{f1} and
\ref{ilcompn}, $\fs_1$ is a skeletal dense homeomorphic embedding
and $\fs_2$ is a skeletal closed homeomorphic embedding. Hence
$\fs$ is a skeletal homeomorphic embedding. Then, by \cite[Theorem
2.11]{Di2}, $f$ is a skeletal homeomorphic embedding.
\sqs

Recall that a continuous mapping $f:X\lra Y$ is said to be {\em
semi-open}\/ (\cite{V2}) if for every point $y\in f(X)$ there
exists a point $x\in f\inv(y)$ such that, for every $U\sbe X$,
$x\in\int_X(U)$ implies that $y\in\int_{f(X)}(f(U))$.

\smallskip

The following assertion is a slight generalization of
\cite[Theorem 6]{F2}.

\begin{theorem}\label{lcaemb}
A morphism $\p\in\SKAL((A,\rho,\BBBB),(B,\eta,\BBBB\ap))$ is an
LCA-embedd\-ing iff the map $f=\Psi^a(\p)$ is a quasi-open
semi-open perfect surjection.
\end{theorem}

\doc Set $X=\Psi^a(A,\rho,\BBBB)$, $Y=\Psi^a(B,\eta,\BBBB\ap)$,
 $C=C_\rho$ and $C\ap=C_\eta$ (see \cite[1.16]{Di2} for
the notations). Then, by \cite[1.16]{Di2}, $(A,C)$ and $(B,C\ap)$
are CNCA's. By the proof of Theorem \cite[2.1]{Di2},
$\Psi^a(A,C)=\a X=X\cup\{\s^A_\infty\}$ and $\Psi^a(B,C\ap)=\a
Y=Y\cup\{\s^B_\infty\}$.

 Let $\p$ be an LCA-embedding, i.e. $\p:A\lra B$ is a Boolean
embedding such that, for any $a,b\in A$, $a\rho b$ iff
$\p(a)\eta\p(b)$, and $a\in\BBBB$ iff $\p(a)\in\BBBB\ap$; hence
$\p$ satisfies condition (L3) from \cite{Di2}. Then, by Theorem
\ref{slcompn}, $f$ is a perfect quasi-open surjection. It remains
to show that $f$ is semi-open.
 Denote by $\p_c$ the map $\p$ regarded as a function from $(A,C)$ to
 $(B,C\ap)$.
  By \cite[(55)]{Di2}, $\p_c$ satisfies
condition (F1)  from \cite[2.12]{Di2}. We will show that $\p_c$ is
an NCA-embedding. Indeed, for any $a,b\in A$, we have that $aCb$
iff $a\rho b$ or $a,b\nin\BBBB$; since $\p$ is an LCA-embedding,
we obtain that $aCb$ iff $\p_c(a)C\ap\p_c(b)$. So, $\p_c$ is an
NCA-embedding and an $\SAC$-morphism. Then, by Theorem 6 of
Fedorchuk's paper \cite{F2}, $f_c=\Psi^a(\p_c):\a Y\lra\a X$ is a
semi-open map. If $1_A\nin\BBBB$ and $1_B\nin\BBBB\ap$ then
$f_c\inv(\s^A_\infty)=\{\s^B_\infty\}$ (see the proof of
\cite[Theorem 2.15]{Di2}) and since $f=f_c|Y$, we obtain that $f$
is semi-open.  Further, if $1_A\in\BBBB$  and $1_B\in\BBBB\ap$
then the fact that $f$ is semi-open is obvious. Since only these
two cases are possible in the given situation, we have proved that
$f$ is a perfect quasi-open semi-open surjection.

Conversely, let $f$ be a perfect quasi-open semi-open surjection.
Then, by Theorem \ref{slcompn}, $\p$ is an injective
$\SAL$-morphism. Hence $\pl\circ\p=id_A$. Thus, if
$\p(a)\in\BBBB\ap$ then, by (L2), $a=\pl(\p(a))\in\BBBB$. Using
(L3), we obtain that $a\in\BBBB$ iff $\p(a)\in\BBBB\ap$. Since
(L1) takes place, it remains only to prove that $a\rho b$ implies
$\p(a)\eta\p(b)$, for all $a,b\in A$. Set $\ps=\Psi^t(f)$. Let
$F,G\in RC(X)$,  $F\cap G\nes$ and $x\in F\cap G$. Set $U=\int(F)$
and $V=\int(G)$. Then $x\in \cl(U)\cap\cl(V)$. Since $f$ is a
semi-open surjection, there exists $y\in f\inv(x)$ such that, for
every $W\sbe Y$, $y\in\int_Y(W)$ implies that $x\in\int_X(f(W))$.
We will show that $y\in \cl(f\inv(U))\cap\cl(f\inv(V))$. Indeed,
suppose that $y\nin \cl(f\inv(U))$. Then there exists an open
neighborhood $Oy$ of $y$ such that $Oy\cap f\inv(U)=\ems$. Thus
$f(Oy)\cap U=\ems$. Since $x\in \cl(U)$ and $x\in \int(f(Oy))$, we
obtain a contradiction. Hence $y\in \cl(f\inv(U))$. Analogously we
can show that $y\in \cl(f\inv(V))$. Therefore, $y\in
\cl(f\inv(U))\cap\cl(f\inv(V))=\ps(F)\cap\ps(G)$. Now, using
\cite[Theorem 2.15]{Di2}, we get that $a\rho b$ implies
$\p(a)\eta\p(b)$.
\sqs

\begin{defi}\label{lideal}
\rm Let $(A,\rho,\BBBB)$ be a CLCA. An ideal $I$ of $A$ is called
a {\em $\d$-ideal} if $I\sbe \BBBB$ and for any $a\in I$ there
exists $b\in I$ such that $a\ll_\rho b$. If $I_1$ and $I_2$ are
two $\d$-ideals of $(A,\rho,\BBBB)$ then we put $I_1\le I_2$ iff
$I_1\sbe I_2$. We will denote by $(I(A,\rho,\BBBB),\le)$ the poset
of all $\d$-ideals of $(A,\rho,\BBBB)$.
\end{defi}

\begin{fact}\label{dideal}
Let $(A,\rho,\BBBB)$ be a CLCA. Then, for every $a\in A$, the set
$\{b\in\BBBB\st b\ll_\rho a\}$ is a $\d$-ideal. Such $\d$-ideals
will be called\/ {\em principal $\d$-ideals}.
\end{fact}

\doc The proof is obvious.
\sqs

Recall that a {\em frame} is a complete lattice $L$ satisfying the
infinite distributive law $a\we\bigvee S=\bigvee\{a\we s\st s\in
S\}$, for every $a\in L$ and every $S\sbe L$.

\begin{fact}\label{frlid}
Let $(A,\rho,\BBBB)$ be a CLCA. Then the poset
$(I(A,\rho,\BBBB),\le)$ of all $\d$-ideals of $(A,\rho,\BBBB)$
(see \ref{lideal}) is a frame.
\end{fact}

\doc It is well known that the set $Idl(A)$ of all ideals of a
distributive lattice forms a frame under the inclusion ordering
(see, e.g., \cite{J2}). It is easy to see that the join in
$(Idl(A),\sbe)$ of a family of $\d$-ideals  is a $\d$-ideal and
hence it is the join of this family in $(I(A,\rho,\BBBB),\le)$.
The meet in $(Idl(A),\sbe)$ of a finite family of $\d$-ideals is
also a $\d$-ideal and hence it is the meet of this family in
$(I(A,\rho,\BBBB),\le)$. Therefore, $(I(A,\rho,\BBBB),\le)$ is a
frame. Note that the meet of an infinite family of $\d$-ideals in
$(I(A,\rho,\BBBB),\le)$ is not obliged to coincide with the meet
of the same family in $(Idl(A),\sbe)$.
\sqs

\begin{theorem}\label{opensetsfr}
Let $(A,\rho,\BBBB)$ be a CLCA, $Y=\Psi^a(A,\rho,\BBBB)$ and
$\OO(Y)$ be the frame of all open subsets of\/ $Y$. Then there
exists a frame isomorphism
$$\iota:(I(A,\rho,\BBBB),\le)\lra (\OO(Y),\sbe),$$
where $(I(A,\rho,\BBBB),\le)$ is the frame of all $\d$-ideals of
$(A,\rho,\BBBB)$. The  isomorphism $\iota$ sends the set
$PI(A,\rho,\BBBB)$ of all principal $\d$-ideals of
$(A,\rho,\BBBB)$ onto the set $RO(Y)$ of all regular open subsets
of\/ $Y$.
\end{theorem}

\doc  Let $I$ be a $\d$-ideal. Put $\iota(I)=\bigcup\{\l^g_A(a)\st a\in
I\}$. Then $\iota(I)$ is an open subset of $Y$. Indeed, for every
$a\in I$ there exists $b\in I$ such that $a\ll b$. Then
$\l^g_A(a)\sbe\int_Y(\l^g_A(b))\sbe \l^g_A(b)\sbe \iota(I)$. Hence
$\iota(I)$ is an open subset of $Y$. Therefore $\iota$ is a
function from $I(A,\rho,\BBBB)$ to $\OO(Y)$. Let $U\in\OO(Y)$. Set
$\BBBB_U=\{b\in\BBBB\st\l^g_A(b)\sbe U\}$. Then, as it is easy to
see, $\BBBB_U$ is a $\d$-ideal of $(A,\rho,\BBBB)$. Since $Y$ is a
locally compact Hausdorff space, $\iota(\BBBB_U)=U$. Hence,
$\iota$ is a surjection. We will show that $\iota$ is an injection
as well. Indeed, let $I_1,I_2\in I(A,\rho,\BBBB)$ and
$\iota(I_1)=\iota(I_2)$. Set $\iota(I_1)=W$ and put
$\BBBB_W=\{b\in\BBBB\st\l^g_A(b)\sbe W\}$. Then, obviously,
$I_1\sbe\BBBB_W$. Further, if $b\in\BBBB_W$ then $\l^g_A(b)\sbe W$
and $\l^g_A(b)$ is compact. Since $I_1$ is a $\d$-ideal,
$\OM=\{\int(\l^g_A(a))\st a\in I_1\}$ is an open cover of $W$ and,
hence, of $\l^g_A(b)$. Thus there exists a finite subcover
$\{\l^g_A(a_1),\ldots,\l^g_A(a_k)\}$ of $\OM$. Therefore,
$\l^g_A(b)\sbe\bigcup\{\l^g_A(a_i)\st
i=1,\ldots,k\}=\l^g_A(\bigvee\{a_i\st i=1,\ldots,k\})$. This
implies that $b\le \bigvee\{a_i\st i=1,\ldots,k\}$ and hence $b\in
I_1$. So, we have proved that $I_1=\BBBB_W$. Analogously we can
show that $I_2=\BBBB_W$. Thus $I_1=I_2$. Therefore, $\iota$ is a
bijection. It is obvious that if $I_1,I_2\in I(A,\rho,\BBBB)$ and
$I_1\le I_2$ then $\iota(I_1)\sbe \iota(I_2)$. Conversely, if
$\iota(I_1)\sbe \iota(I_2)$ then $I_1\le I_2$. Indeed, if
$\iota(I_i)=W_i,\ i=1,2$, then, as we have already seen,
$I_i=\BBBB_{W_i},\ i=1,2$; since $W_1\sbe W_2$ implies that
$\BBBB_{W_1}\sbe\BBBB_{W_2}$, we get that $I_1\le I_2$. So,
$\iota:(I(A,\rho,\BBBB),\le)\lra(\OO(Y),\sbe)$ is an isomorphism
of posets. This implies that $\iota$ is also a frame isomorphism.

Let $U$ be a regular open subset of $Y$. Set $F=Y\stm U$. Then
there exists $a\in A$ such that $F=\l^g_A(a)$. Put
$\BBBB_U=\{b\in\BBBB\st \l^g_A(b)\sbe U\}$. Then, as we have
already seen, $\BBBB_U$ is a $\d$-ideal and $\iota(\BBBB_U)=U$.
Since $F\in RC(Y)$, we have that $\BBBB_U=\{b\in\BBBB\st b\ll_\rho
a^*\}$. Hence $\BBBB_U$ is a principal $\d$-ideal.

Conversely, if $I$ is a principal $\d$-ideal then $U=\iota(I)$ is
a regular open set in $Y$. Indeed, let $a\in A$ and $I=\{b\in
\BBBB\st b\ll_\rho a\}$. It is enough to prove that $Y\stm
U=\l^g_A(a^*)$. If $b\in I$ then $b(-\rho) a^*$ and hence
$\l^g_A(b)\cap\l^g_A(a^*)=\ems$. Thus $U\sbe Y\stm \l^g_A(a^*)$.
If $\s\in Y\stm\l^g_A(a^*)$ then, by \cite[(24)]{Di2}, there
exists $b\in \BBBB$ such that $\s\in\l^g_A(b)\sbe Y\stm
\l^g_A(a^*)$. Since, by \cite[(11)]{Di2} and \cite[(42)]{Di2},
$Y\stm \l^g_A(a^*)=\int_Y(\l^g_A(a))$, we get that $b\ll_\rho a$.
Therefore $b\in I$ and hence $\s\in U$. This means that
$Y\stm\l^g_A(a^*)\sbe U$.
\sqs

Recall that if $A$ is a Boolean algebra and $a\in A$ then the set
$\downarrow a$ endowed with the same meets and joins as in $A$ and
with complements $\neg b$ defined by the formula $\neg b=b^*\we
a$, for every $b\le a$, is a Boolean algebra; it is denoted by
$A|a$. If $J=\downarrow a^*$ then $A|a$ is isomorphic to the
factor algebra $A/J$; the isomorphism $h:A|a\lra A/J$ is the
following: $h(b)=[b]$, for every $b\le a$ (see, e.g., \cite{Si2}).

With the next theorem we show how one can build the CLCAs
corresponding to the open subsets of a locally compact Hausdorff
space $Y$ from the CLCA $\Psi^t(Y)$.

\begin{theorem}\label{conopen}
Let $(A,\rho,\BBBB)$ be a CLCA, $Y=\Psi^a(A,\rho,\BBBB)$ and $I$
be a $\d$-ideal of $(A,\rho,\BBBB)$. Set $a_I=\bigvee I$ and
$B=A|a_I$. For every $a,b\in B$, set $a\eta b$ iff there exist
$c\in I$ and $\s\in Y$ such that $a,b,c\in\s$. Then $(B,\eta,I)$
is a CLCA. Let $\p:A\lra B$ be the natural epimorphism (i.e.
$\p(a)=a\we a_I$, for every $a\in A$). If $X=\Psi^a(B,\eta,I)$ and
$f=\Psi^a(\p)$ then $f:X\lra Y$ is an open injection and
$f(X)=\iota(I)$ (see \ref{opensetsfr} for $\iota$). (Hence, $X$ is
homeomorphic to $\iota(I)$.)
\end{theorem}

\doc Obviously, $B$ is a  complete Boolean algebra and
$\p$ is a surjective complete Boolean homomorphism.

Set $U=\iota(I)$ (i.e. $U=\bigcup\{\l_A^g(b)\st b\in I\}$). Then
$U$ is open in $Y$ (see \ref{opensetsfr}) and
$\cl(U)=\l_A^g(a_I)$. Since $I$ is a $\d$-ideal,
$\{\int(\l_A^g(b))\st b\in I\}$ is an open cover of $U$.

If $I=\{0\}$ then $U=\ems$, $a_I=0$, $B=\{0\}$ and $X=\ems$;
hence, in this case the assertion of the theorem is true. Thus,
let us assume that $I\neq \{0\}$.

We will first check that $(B,\eta,I)$ is a CLCA, i.e. that
conditions  (C1)-(C4)  and (BC1)-(BC3) from \cite{Di2} are
fulfilled.

Let $b\in B\stm \{0\}$. Set $V=\int(\l_A^g(b))$. Since $b\le a_I$,
we obtain that $\ems\neq V\sbe\l_A^g(b)\sbe\l_A^g(a_I)=\cl(U)$.
Hence $V\cap U\nes$. Since $Y$ is locally compact, there exists
$c\in \BBBB\stm \{0\}$ such that $\l_A^g(c)\sbe V\cap U$. Thus
$\l_A^g(c)\sbe U$ and, hence, there exists $d\in I$ such that
$\l_A^g(c)\sbe \l_A^g(d)$ (we use the fact that $\l_A^g(c)$ is
compact and $I$ is closed under finite joins). We get that $c\le
d$. Therefore, $c\in I$. There exists $\s\in Y$ such that
$c\in\s$. Since $c\le b$, we obtain that $b\in\s$. Hence, $b\eta
b$. So, the axiom (C1) is fulfilled. Note that $c\ll_\rho b$
(because $\l_A^g(c)\sbe V=\int(\l_A^g(b))$). Thus, $c\ll_\eta b$.
Therefore, the axiom (BC3) is checked as well. Clearly, the axioms
(C2) and (C3) are satisfied. Using condition (K2) from \cite{Di2},
we obtain that the axiom (C4) is also fulfilled.

 Let $a\in I$, $c\in B$ and $a\ll_\eta c$.  Then, for every $\s\in Y$,
 we have that either
$a\nin\s$ or $(c^*\we a_I)\nin\s$. Since $a\in\BBBB$, we get that
$a(-C_\rho)(c^*\we a_I)$ (see   \cite[Lemma 20]{VDDB2}). Using
again the fact that $a\in\BBBB$, we obtain that $a\ll_\rho (c\vee
a_I^*)$. Then there exists $b\in\BBBB$ such that $a\ll_\rho
b\ll_\rho (c\vee a_I^*)$. Since $I$ is a $\d$-ideal, there exists
$d\in I$ such that $a\ll_\rho d$. Then $a\ll_\rho b\we d$ and
$b\we d\in I$. Thus  $a\ll_\eta b\we d\ll_\eta c$. Therefore, the
axiom (BC1) is checked.

Let us note that
\begin{equation}\label{equ}
\mbox{if } \s \mbox{ is a cluster in } (A,C_\rho),\ a\in\s \mbox{
and } b^*\nin\s \mbox{ then } a\we b\in\s.
\end{equation}
Indeed, there exists an ultrafilter $u$ such that $a\in u\sbe \s$.
Then $b^*\nin u$ and hence $b\in u$. Therefore $a\we b\in u$. Thus
$a\we b\in\s$. So, (\ref{equ}) is proved.

 Let $a,b\in B$ and $a\eta b$. Then there exist $c\in I$ and $\s\in Y$ such
that $a,b,c\in\s$.  Since $I$ is a $\d$-ideal, there exists $d\in
I$ such that $c\ll_\rho d$. Then $d^*\nin\s$.  Therefore, by
(\ref{equ}), $b\we d\in\s$. Since $a,b\we d\in\s$, we obtain that
$a\eta(b\we d)$. Hence, the axiom (BC2) is fulfilled.

So, we have proved that $(B,\eta,\BBBB\ap)$ is a CLCA.

We will show that $\p$ satisfies axioms (L1), (L2) and (LO) from
\cite{Di2}. Note first that, for every $a\in B$,
\begin{equation}\label{eq1}
\pl(a)=a.
\end{equation}
This observation  shows that $\p$ satisfies condition (L2). For
checking condition (EL1) from \cite{Di2} (which is equivalent to
the condition (L1)), let $a,b\in B$ and $a\eta b$. Then there
exist $c\in I$ and $\s\in Y$ such that $a,b,c\in\s$. Since $I$ is
a $\d$-ideal, there exists $d\in I$ such that $c\ll_\rho d$. Then
$d^*\nin\s$ and (\ref{equ}) implies that $a\we d,b\we d\in\s$.
Since $a\we d,b\we d\in\BBBB$, we get that $(a\we d)\rho (b\we
d)$. Therefore, $a\rho b$. So, $\p$ satisfies condition (EL1). Let
us prove that the axiom (LO) is fulfilled as well. Let $a\in A$,
$b\in I$ and $a\rho b$. Then $a C_\rho b$ and since $b\in \BBBB$,
Lemma 20 from \cite{VDDB2} implies that there exists $\s\in Y$
such that $a,b\in\s$. There exists $d\in I$ such that $b\ll_\rho
d$ (because $I$ is a $\d$-ideal). Then $b\ll_{C_\rho} a_I$. Thus
$a_I^*\nin\s$. Now, (\ref{equ}) implies that $a\we a_I\in\s$.
Hence $\p(a)\eta b$. Therefore,  condition (LO) is checked. So,
$\p$ satisfies axioms (L1), (L2) and (LO). Then, by Theorem
\ref{ilcompo}, $f:X\lra Y$ is an open injection and hence $f$ is a
homeomorphism between $X$ and $f(X)$. Let us show that $f(X)=U$.

Set $\ps=\Psi^t(f)$. Then $\ps:RC(Y)\lra RC(X)$, $\ps(F)=f\inv(F)$
for every $F\in RC(Y)$ and $\l_B^g\circ\p=\ps\circ\l_A^g$ (see
\cite[Theorem 2.17]{Di2}). Hence, for every $b\in I$,
$\l_B^g(b)=f\inv(\l_A^g(b))\sbe f\inv(U)$. Since
$X=\bigcup\{\l_B^g(b)\st b\in I\}$, we obtain that $f(X)\sbe U$.
For showing that $f(X)\spe U$, let $\s\in U$. Then, by the
definition of $U$, there exists $b\in I$ such that
$\s\in\l_A^g(b)\sbe U$. Hence $b\in\s$ and $a_I^*\nin\s$. There
exists $u\in Ult(A)$ such that $b\in u\sbe \s$ and $\s=\s_u$.
Since $a_I^*\nin\s$, we obtain that $a_I\in u$. Thus, for every
$c\in u$, $c\we a_I\neq 0$. Therefore, $\p(c)\neq 0$, for every
$c\in u$. This implies that $\p(u)$ is a basis of a bounded filter
in $B$. Then there exists $v\in Ult(B)$ such that $\p(u)\sbe v$.
Hence $u\sbe\p\inv(v)$ and thus $u=\p\inv(v)$. Since $\s_v$ is a
bounded cluster, we get that $\s_v\in X$. Further,
$f(\s_v)=\s_{\p\inv(v)}=\s_u=\s$. Hence $f(X)=U$.
\sqs

We will now show how one can build the CLCAs corresponding to the
regular closed subsets of a locally compact Hausdorff space $Y$
from the CLCA $\Psi^t(Y)$.

\begin{theorem}\label{conregclo}
Let $(A,\rho,\BBBB)$ be a CLCA, $Y=\Psi^a(A,\rho,\BBBB)$, $a_0\in
A$ and $F=\l_A^g(a_0)$. Set $B=A|a_0$ and let $\p:A\lra B$ be the
natural epimorphism (i.e. $\p(a)=a\we a_0$, for every $a\in A$).
Put $\BBBB\ap=\p(\BBBB)$ and let, for every $a,b\in B$, $a\eta b$
iff $a \rho b$. Then $(B,\eta,\BBBB\ap)$ is a CLCA. If
$X=\Psi^a(B,\eta,\BBBB\ap)$ and $f=\Psi^a(\p)$ then $f:X\lra Y$ is
a closed quasi-open injection and $f(X)=F$. (Hence, $X$ is
homeomorphic to $F$.)
\end{theorem}

\doc We have that $B$
is a complete Boolean algebra, $\p$ is a complete Boolean
homomorphism and $\pl(a)=a$, for every $a\in B$. Set
$\psi=\l_A^g\circ\pl$. We will show that
$\psi\ap=\psi\upharpoonright B$ is a CLCA-isomorphism between
$(B,\eta,\BBBB\ap)$ and $(RC(F),\rho_F,CR(F))$. Since $F\in
RC(Y)$, we have, as it is well known, that $RC(F)\sbe RC(Y)$ and
$RC(F)=\{G\we F\st G\in RC(Y)\}$; moreover, $RC(F)=RC(Y)|F$.
 Hence $\psi\ap:B\lra RC(F)$ is a Boolean isomorphism.  For any
$a,b\in B$, we have that $a\eta b\iff \pl(a)\rho\pl(b)\iff
\l_A^g(\pl(a))\cap\l_A^g(\pl(b))\nes\iff \psi(a)\rho_F\psi(b)$.
Finally,  for any $a\in B$, we have that $a\in\BBBB\ap\iff
\pl(a)\in\BBBB\iff \l_A^g(\pl(a))$ is compact $\iff \psi(a)\in
CR(F)$. Therefore, $(B,\eta,\BBBB\ap)$ is a CLCA; it is isomorphic
to the CLCA $(RC(F),\rho_F,CR(F))$. For showing that $f:X\lra Y$
is a homeomorphic embedding and $f(X)=F$, note that $\p$ satisfies
conditions (L1)-(L3) from \cite{Di2} and condition (LS), and
hence, by Theorem \ref{ilcompn}, $f$ is a quasi-open perfect
injection, i.e. $f$ is a homeomorphic embedding. From
\cite[(45)]{Di2} we get that, for every $b\in\BBBB\ap$,
$f(\l_B^g(b))=\l_A^g(\pl(b))=\l_A^g(b)\sbe F$. Since
$X=\bigcup\{\l_B^g(b)\st b\in\BBBB\ap\}$, we obtain that $f(X)\sbe
F$. Let $y\in\int_Y(F)$. Then there exists $b\in\BBBB$ such that
$y\in\int(\l_A^g(b))\sbe \l_A^g(b)\sbe\int_Y(F)$. Hence
$b\in\BBBB\ap$. Using again \cite[(45)]{Di2}, we get that $y\in
f(\l_B^g(b))$, i.e. $y\in f(X)$. Thus $\int_Y(F)\sbe f(X)$. Since
$f(X)$ is closed in $Y$, we conclude that $f(X)\spe F$. Therefore,
$f(X)=F$.
\sqs

\section{Some applications}

We start with a proposition which has a straightforward proof.

\begin{pro}\label{isopo}
Let $(A,\rho,\BBBB)$ be a CLCA, $Y=\Psi^a(A,\rho,\BBBB)$ and $a\in
A$. Then $a$ is an atom of $A$ iff  $\l_A^g(a)$ is an isolated
point of the space $Y$.
\end{pro}

\begin{cor}\label{isopocor}
Let $(A,\rho,\BBBB)$ be a CLCA. Then\/ $\BBBB$ contains all finite
sums of the atoms of $A$.
\end{cor}

\begin{pro}\label{disspa}
Let $(A,\rho,\BBBB)$ be a CLCA and $Y=\Psi^a(A,\rho,\BBBB)$. Then
$Y$ is a discrete space  iff\/ $\BBBB$ is the set of all finite
sums of the atoms of $A$.
\end{pro}

\doc It follows easily from \ref{isopo} and the fact that
$Y=\bigcup\{\l_A^g(a)\st a\in\BBBB\}$. \sqs

The next proposition has an easy proof which will be omitted.

\begin{pro}\label{extrspa}
Let $(A,\rho,\BBBB)$ be a CLCA and $Y=\Psi^a(A,\rho,\BBBB)$. Then
$Y$ is an extremally disconnected  space  iff $a\llx a$, for every
$a\in A$.
\end{pro}

\begin{pro}\label{denseiso}
Let $(A,\rho,\BBBB)$ be a CLCA and $Y=\Psi^a(A,\rho,\BBBB)$. Then
the set of all isolated points of $Y$ is dense in $Y$ iff $A$ is
an atomic Boolean algebra.
\end{pro}

\doc ($\Rightarrow$) Let $X$ be the set of all isolated points of
$Y$ and let $\cl_Y(X)=Y$. Then, by \ref{isombool}, the Boolean
algebras $RC(X)$ and $RC(Y)$ are isomorphic.  From \cite[Theorem
2.11]{Di2} we get that $A$ is isomorphic to $RC(Y)$ and, hence, to
$RC(X)$. Therefore, $A$ is an atomic Boolean algebra.

($\Leftarrow$) Set $X=\{\l_A^g(a)\st a\in Atoms(A)\}$. Then, by
Proposition \ref{isopo},  $X$ consists of isolated points of $Y$.
We need only to show that $X$ is dense in $Y$. Let $b\in\BBBB$ and
$b\neq 0$. Then there exists $a\in Atoms(A)$ such that $a\le b$.
Thus $\l_A^g(a)\in\l_A^g(b)$, i.e. $X\cap\l_A^g(b)\nes$. Since $Y$
is a regular space, \cite[(24)]{Di2} implies that $X$ is dense in
$Y$. \sqs

\begin{nota}\label{locs}
\rm Let $X$ be a Tychonoff space. We will denote by $\LL(X)$ the
set of all, up to equivalence, locally compact Hausdorff
extensions of $X$ (recall that two (locally compact Hausdorff)
extensions $(Y_1,f_1)$ and $(Y_2,f_2)$ of $X$ are said to be {\em
equivalent}\/ iff there  exists a homeomorphism $h:Y_1\lra Y_2$
such that $h\circ f_1=f_2$). We will regard two orders on it. Let
$[(Y_i,f_i)]\in\LL(X)$, where $i=1,2$. We set $[(Y_1,f_1)]\le
[(Y_2,f_2)]$ (respectively, $[(Y_1,f_1)]\le_s [(Y_2,f_2)]$) iff
there exists a continuous (resp., continuous surjective) mapping
$h:Y_2\lra Y_1$ such that $f_1=h\circ f_2$.
\end{nota}

\begin{defi}\label{locsd}
\rm Let $X$ be a locally compact Hausdorff space and
$\Psi^t(X)=(A,\rho,\BBBB)$. We will denote by $\LL_a(X)$ the set
of all LCAs of the form $(A,\rho_1,\BBBB_1)$ which satisfy the
following conditions:

\noindent(LA1) $\rho\sbe\rho_1$;

\noindent(LA2) $\BBBB\sbe\BBBB_1$;

\noindent(LA3) for every $a\in A$ and every $b\in\BBBB$, $b\rho_1
a$ implies $b\rho a$.

We will define two orders on the set $\LL_a(X)$. If
$(A,\rho_i,\BBBB_i)\in\LL_a(X)$, where $i=1,2$, we set
$(A,\rho_1,\BBBB_1)\preceq(A,\rho_2,\BBBB_2)$ (respectively,
$(A,\rho_1,\BBBB_1)\preceq_s (A,\rho_2,\BBBB_2)$)     iff
$\rho_2\sbe\rho_1$ and $\BBBB_2\sbe\BBBB_1$ (and, respectively, in
addition, for every bounded ultrafilter $u$ in
$(A,\rho_1,\BBBB_1)$ there exists $b\in\BBBB_2$ such that $b\rho_1
u$).
\end{defi}

\begin{theorem}\label{locset}
Let $(X,\tau)$ be a locally compact Hausdorff space and let
$\Psi^t(X,\tau)=(A,\rho,\BBBB)$. Then there exists an isomorphism
$\mu$ (respectively, $\mu_s$) between the ordered sets
$(\LL(X),\le)$ and  $(\LL_a(X),\preceq)$ (respectively,
$(\LL(X),\le_s)$ and  $(\LL_a(X),\preceq_s)$).
\end{theorem}

\doc If $X$ is compact then everything is clear. Thus, let $X$ be a
non-compact space.

Let $[(Y,f)]\in\LL(X)$. Then $f:X\lra Y$ is a dense homeomorphic
embedding and hence it is an open injection . Set
$(A\ap,\rho\ap,\BBBB\ap)=\Psi^t(Y)$ and $\p=\Psi^t(f)$. Then, by
Theorem \ref{f1}, $\p:(A\ap,\rho\ap,\BBBB\ap)\lra (A,\rho,\BBBB)$
is a Boolean isomorphism and a $\SOAL$-morphism. For every $a,b\in
A$, set $a\rho_f b$ iff $\p\inv(a)\rho\ap \p\inv(b)$, and set
$a\in\BBBB_f$ iff $\p\inv(a)\in\BBBB\ap$. For every $a\in A$, set
$\psi(a)=\p\inv(a)$. Then $\psi:(A,\rho_f,\BBBB_f)\lra (A\ap,
\rho\ap,\BBBB\ap)$ is an LCA-isomorphism. It is easy to see that
$(A,\rho_f,\BBBB_f)\in\LL_a(X)$. Set
$\mu([(Y,f)])=(A,\rho_f,\BBBB_f)$. Then it is not difficult to
show that $\mu$ is a well-defined order preserving map between the
ordered sets $(\LL(X),\le)$ and  $(\LL_a(X),\preceq)$.

Let $(A,\rho_1,\BBBB_1)\in\LL_a(X)$. Then the identity map
$i:(A,\rho_1,\BBBB_1)\lra(A,\rho,\BBBB)$ is a Boolean isomorphism
and a $\SOAL$-morphism. Set $X\ap=\Psi^a(A,\rho,\BBBB)$,
$Y=\Psi^a(A,\rho_1,\BBBB_1)$ and $f\ap=\Psi^a(i)$. Then, by
Theorem \ref{f1}, $f\ap: X\ap\lra Y$ is a dense homeomorphic
embedding. Set $f=f\ap\circ t_X$ (see \cite[(26)]{Di2} for the
definition of the homeomorphism $t_X:X\lra X\ap$). Then
$[(Y,f)]\in\LL(X)$ and we put $\mu\ap(A,\rho_1,\BBBB_1)=[(Y,f)]$.
It is easy to show that
$\mu\ap:(\LL_a(X),\preceq)\lra(\LL(X),\le)$ is an order preserving
map.

Using \cite[Theorem 2.11]{Di2}, it is not difficult to prove that
the compositions $\mu\circ\mu\ap$ and $\mu\ap\circ\mu$ are
identities.

Finally, we will show that the same map $\mu$, which will be now
denoted  by $\mu_s$, is an isomorphism between the ordered sets
$(\LL(X),\le_s)$ and  $(\LL_a(X),\preceq_s)$. This can be proved
easily using Theorem \ref{slcompnk} and \cite[Proposition
3.2]{DV2} (the last proposition says that if $(B,\eta)$ is a CA
and $F_1$, $F_2$ are two filters in $B$ such that $F_1\eta F_2$
then there exist ultrafilters $u_1,u_2$ in $B$ such that $F_i\sbe
u_i$, where $i=1,2$, and $u_1\eta u_2$).

This completes the proof of our theorem. \sqs

Recall that if $X$ is a set and $\PP(X)$ is the power set of $X$
ordered by the inclusion, then a triple $(X,\rho,\BBBB)$ is called
a {\em local proximity space} (see \cite{LE2}) if $(\PP(X),\rho)$
is a CA, $\BBBB$ is an ideal (possibly non proper) of $\PP(X)$ and
the axioms (BC1),(BC2) from \cite[Definition 1.14]{Di2} are
fulfilled. A local proximity space $(X,\rho,\BBBB)$ is said to be
{\em separated} if $\rho$ is the identity relation on singletons.

\begin{rem}\label{remd}
\rm In this remark we will use the notations from Theorem
\ref{locset}. Let $(A,\rho_1,\BBBB _1)\in\LL_a(X)$ and
$\mu\inv(A,\rho_1,\BBBB_1)=[(Y,f)](\in\LL(X))$. Let
 $(X,\rho_Y,\BBBB_Y)$ be the
local proximity space  induced by $(Y,f)$ (i.e., for every
$M,N\sbe X$, $M\rho_Y N \iff \cl_Y(f(M))\cap\cl_Y(f(N))\nes$ and
$M\in\BBBB_Y$ iff $\cl_Y(f(M))$ is compact). Then
$\BBBB_1=\BBBB_Y\cap RC(X)$ and, for every $F,G\in RC(X)$,
$F\rho_Y G\iff F\rho_1 G$. Indeed, using \cite[Theorem 2.11 and
(31)]{Di2} and the fact that $f\ap$ is an open mapping, we obtain
that, for every $F\in RC(X)$,
$\l^g_{(A,\rho,\BBBB)}(F)=(f\ap)\inv(\l^g_{(A,\rho_1,\BBBB_1)}(F))$
and thus, by \ref{isombool},
$\cl_Y(f\ap(\l^g_{(A,\rho,\BBBB)}(F)))=
\cl_Y(f\ap((f\ap)\inv(\l^g_{(A,\rho_1,\BBBB_1)}(F))))
=\l^g_{(A,\rho_1,\BBBB_1)}(F)$. As it follows from \cite[Theorem
2.11]{Di2}, $\Psi^t(t_X)\circ \l^g_{(A,\rho,\BBBB)}=id_A$. Hence,
using \cite[(31)]{Di2}, we obtain that, for every $F\in RC(X)$,
$F=(\Psi^t(t_X))(\l^g_{(A,\rho,\BBBB)}(F))=
\cl_{X\ap}(\int_{X\ap}(t_X\inv(\l^g_{(A,\rho,\BBBB)}(F))))$ and
thus $F=t_X\inv(\l^g_{(A,\rho,\BBBB)}(F))$, i.e.
$t_X(F)=\l^g_{(A,\rho,\BBBB)}(F)$. Therefore,
$\cl_Y(f(F))=\l^g_{(A,\rho_1,\BBBB_1)}(F)$ for every $F\in RC(X)$.
Now, all follows from \cite[(III) and (IV) in the proof of Theorem
2.1]{Di2}.
\end{rem}

\begin{nota}\label{algcomp}
\rm If $(A,\rho,\BBBB)$ is a CLCA then we will write
$\rho\sbe_{\BBBB} C$ provided that $C$ is a normal contact
relation on $A$ satisfying the following conditions:

\noindent(RC1) $\rho\sbe C$, and

\noindent(RC2) for every $a\in A$ and every $b\in\BBBB$, $aCb$
implies $a\rho b$.

\smallskip

\noindent If $\rho\sbe_{\BBBB} C_1$ and $\rho\sbe_{\BBBB} C_2$
then we will write $C_1\preceq_c C_2$ iff $C_2\sbe C_1$.
\end{nota}

\begin{cor}\label{compset}
Let $(X,\tau)$ be a locally compact Hausdorff space and let
$\Psi^t(X,\tau)=(A,\rho,\BBBB)$. Then there exists an isomorphism
$\mu_c$ between the ordered set $(\KK(X),\le)$ of all, up to
equivalence, Hausdorff compactifications of $(X,\tau)$ and the
ordered set $(\KK_a(X),\preceq_c)$ of all normal contact relations
$C$ on $A$ such that $\rho\sbe_{\BBBB} C$ (see \ref{algcomp} for
the notations).
\end{cor}

\doc It follows immediately from Theorem \ref{locset}.
 \sqs

\begin{pro}\label{alphabeta}
Let $(X,\tau)$ be a locally compact non-compact Hausdorff space
and let $\Psi^t(X,\tau)=(A,\rho,\BBBB)$. Then $C_\rho$ is the
smallest element of the ordered set $(\KK_a(X),\preceq_c)$ (see
\cite[Lemma 1.16]{Di2} for $C_\rho$); hence, if $(\a X,\a)$ is the
Alexandroff (one-point) compactification  of $X$ then $\mu_c([(\a
X,\a)])=C_\rho$ (see Corollary \ref{compset}  for $\mu_c$).
Further,  the ordered set $(\KK_a(X),\preceq_c)$ has a greatest
element $C_{\b\rho}$; it is defined as follows: for every $a,b\in
A$, $a(-C_{\b\rho}) b$ iff there exists a set $\{c_d\in A\st
d\in\DDDD\}$  such that:

\noindent(1) $a\llx c_d\llx b^*$, for all $d\in\DDDD$, and

\noindent(2) for any two elements $d_1,d_2$ of $\DDDD$, $d_1<d_2$
implies that $c_{d_1}\llx c_{d_2}$.

Hence,  if $(\b X,\b)$ is the Stone-\v{C}ech compactification of
$X$  then $\mu_c([(\b X,\b)])=C_{\b\rho}$.
\end{pro}

\doc Recall that, for every $a,b\in A$, $aC_\rho b$ iff either
$a\rho b$ or $a,b\nin\BBBB$.  Obviously, $C_\rho\in\KK_a(X)$. It
is easy to see that if $C\in\KK_a(X)$ then $C\sbe C_\rho$. Hence,
by Corollary \ref{compset}, $\mu_c([(\a X,\a)])=C_\rho$.

The proof that $C_{\b\rho}\in\KK_a(X)$ is straightforward. Let
$C\in\KK_a(X)$. We will show that $C_{\b\rho}\sbe C$. Let $a,b\in
A$ and $aC_{\b\rho}b$. Suppose that $a(-C)b$. Then $a\ll_C b^*$.
Hence there exists $c_{\frac{1}{2}}\in A$ such that $a\ll_C
c_{\frac{1}{2}}\ll_C b^*$.  Since $\rho\sbe C$, we obtain that
$a\llx c_{\frac{1}{2}}\llx b^*$.   It is clear now that we can
construct a set $\{c_d\in A\st d\in\DDDD\}$ such that $a\llx
c_d\llx b^*$, for all $d\in\DDDD$, and for any two elements
$d_1,d_2$ of $\DDDD$, $d_1<d_2$ implies that $c_{d_1}\llx
c_{d_2}$. Thus, $a(-C_{\b\rho})b$, a contradiction. Therefore,
$aCb$. Now, Corollary \ref{compset} implies that $\mu_c([(\b
X,\b)])=C_{\b\rho}$. \sqs

\begin{rem}\label{rembeta}
\rm The definition of the relation $C_{\b\rho}$ in Proposition
\ref{alphabeta} is given in the language of contact relations. It
is clear that if we use the fact that all happens in a topological
space $X$ then we can  define the relation $C_{\b\rho}$ by setting
for every $a,b\in A$, $a(-C_{\b\rho}) b$ iff $a$ and $b$ are
completely separated.
\end{rem}

\begin{pro}\label{joins}
Let $X$ be a locally compact non-compact Hausdorff space. Let
$\Psi^t(X,\tau)=(A,\rho,\BBBB)$ and let $\{C_m\st m\in M\}$ be a
subset of $\KK_a(X)$ (see \ref{compset} for $\KK_a(X)$). For every
$a,b\in A$, put $a(-C)b$ iff there exists a set $\{c_d\in A\st
d\in\DDDD\}$ such that:

\noindent(1) $a\ll_{C_m} c_d\ll_{C_m} b^*$, for all $d\in\DDDD$
and for each $m\in M$, and

\noindent(2) for any two elements $d_1,d_2$ of $\DDDD$, $d_1<d_2$
implies that $c_{d_1}\ll_{C_m} c_{d_2}$, for every $m\in M$.

\noindent Then $C$ is the supremum in $(\KK_a(X),\preceq_c)$ of
the set $\{C_m\st m\in M\}$.
\end{pro}

\doc The proof is straightforward. \sqs
% Podrobnoto d-vo e na listi, prikrepeni na gxrba na Svitxk 5,5.

\begin{rem}\label{discr}
\rm Note that if $X$ is an infinite discrete space and
$\Psi^t(X)=(A,\rho,\BBBB)$ then $A$ is the power set of $X$
ordered by the inclusion, $\BBBB$ is the set of all finite subsets
of $X$ and $\rho$ is the smallest normal contact relation on $A$.
Hence, in this case, $\rho=C_{\b\rho}$.

Note also that Propositions \ref{joins} and \ref{alphabeta} imply
that $(\KK_a(X),\preceq_c)$ is a complete lattice. Now, using
Corollary \ref{compset}, we obtain a new proof of the well-known
fact that, for every locally compact space $X$, $(\KK(X),\le)$ is
a complete lattice.
\end{rem}

 We are now going  to strengthen the Leader Local Compactification
 Theorem (\cite{LE2}) in the same manner as de Vries (\cite{dV2})
 strengthen Smirnov Compactification Theorem (\cite{Sm2}).

Recall that every local proximity space $(X,\rho,\BBBB)$ induces a
completely regular topology $\tau_{(X,\rho,\BBBB)}$ in $X$ by
defining $\cl(M)=\{x\in X\st x\rho M\}$ for every $M\sbe X$
(\cite{LE2}). If $(X,\tau)$ is a topological space then we say
that $(X,\rho,\BBBB)$ is a {\em local proximity space on}
$(X,\tau)$ if $\tau_{(X,\rho,\BBBB)}=\tau$.

\begin{lm}\label{lead1}
Let $(X,\rho_i,\BBBB_i)$, $i=1,2$, be two separated local
proximity spaces on a Tychonoff space $(X,\tau)$ such that
$\BBBB_1\cap RC(X)=\BBBB_2\cap RC(X)$ and $\rho_1|_{RC(X)}=
\rho_2|_{RC(X)}$ (i.e., for every $F,G\in RC(X)$, $F\rho_1 G\iff
F\rho_2 G$). Then $\rho_1=\rho_2$ and $\BBBB_1=\BBBB_2$.
\end{lm}

\doc Let $(Y_i,f_i)$ be the locally compact extension of $X$
generated by $(X,\rho_i,\BBBB_i)$, where $i=1,2$ (see \cite{LE2}).
Let $B\in\BBBB_1$. Then $\cl_{Y_1}(f_1(B))$ is compact. There
exists an open subset $U$ of $Y_1$ such that
$\cl_{Y_1}(f_1(B))\sbe U$ and $\cl_{Y_1}(U)$ is compact. Let
$F=f_1\inv(\cl_{Y_1}(U))$. Then $F\in RC(X)$ and
$\cl_{Y_1}(f_1(F))=\cl_{Y_1}(U)$. Hence  $F\in\BBBB_1\cap RC(X)$.
Thus $F\in\BBBB_2$. Since $B\sbe F$, we get that $B\in\BBBB_2$.
Therefore, $\BBBB_1\sbe\BBBB_2$. Analogously we obtain that
$\BBBB_2\sbe \BBBB_1$. Thus $\BBBB_1=\BBBB_2$.

Let $M,N\sbe X$ and $M(-\rho_1)N$. Suppose that $M\rho_2 N$. Then
there exist $M\ap,N\ap\in\BBBB_2$ such that $M\ap\sbe M$,
$N\ap\sbe N$ and $M\ap\rho_2 N\ap$. Since $\BBBB_1=\BBBB_2$, we
get that $M\ap,N\ap\in\BBBB_1$. Hence $K_1=\cl_{Y_1}(f_1(M\ap))$
and $K_2=\cl_{Y_1}(f_1(N\ap))$ are disjoint compact subsets of
$Y_1$. Then there exist open subsets $U$ and $V$ of $Y_1$ having
disjoint closures in $Y_1$ and containing, respectively, $K_1$ and
$K_2$. Set $F=f_1\inv(\cl_{Y_1}(U))$ and
$G=f_1\inv(\cl_{Y_1}(V))$. Then $F,G\in RC(X)$, $M\ap\sbe F$,
$N\ap\sbe G$ and $F(-\rho_1)G$. Thus $F(-\rho_2)G$ and hence
$M\ap(-\rho_2)N\ap$, a contradiction. Therefore, $M(-\rho_2)N$.
So, $\rho_2\sbe\rho_1$. Using the symmetry, we obtain that
$\rho_1=\rho_2$. \sqs

\begin{defi}\label{lead2}
\rm  Let $(X,\tau)$ be a Tychonoff space. An LCA
$(RC(X,\tau),\rho,\BBBB)$ is said to be {\em admissible for}
$(X,\tau)$ if it satisfies the following conditions:

\noindent(A1) if $F,G\in RC(X)$ and $F\cap G\nes$ then $F\rho G$
(i.e. $\rho_X\sbe\rho$ (see \cite[1.10]{Di2} for $\rho_X$));

\noindent(A2) if $F\in RC(X)$ and $x\in\int_X(F)$ then there
exists $G\in\BBBB$ such that $x\in\int_X(G)$ and $G\llx F$.

The set of all LCAs which are admissible for $(X,\tau)$ will be
denoted by $\LL_{ad}(X,\tau)$ (or simply by $\LL_{ad}(X)$). If
$(RC(X),\rho_i,\BBBB_i)\in\LL_{ad}(X)$, where $i=1,2$, then we set
$(RC(X),\rho_1,\BBBB_1)\preceq_l (RC(X),\rho_2,\BBBB_2)$ iff
$\rho_2\sbe\rho_1$ and $\BBBB_2\sbe\BBBB_1$.
\end{defi}

\begin{lm}\label{lead3}
Let $(X,\rho,\BBBB)$ be a separated local proximity space. Set
$\tau=\tau_{(X,\rho,\BBBB)}$. Let $\rho\ap=\rho|_{RC(X,\tau)}$ and
$\BBBB\ap=\BBBB\cap RC(X,\tau)$. Then
$(RC(X,\tau),\rho\ap,\BBBB\ap)\in\LL_{ad}(X,\tau)$.
\end{lm}

\doc The fact that $(RC(X,\tau),\rho\ap,\BBBB\ap)$ is an LCA is
proved in \cite[Example 40]{VDDB2}. The rest can be easily
checked. \sqs

\begin{lm}\label{lead4}
Let $(X,\tau)$ be a Tychonoff space and
$(RC(X),\rho\ap,\BBBB\ap)\in\LL_{ad}(X)$. Then there exists a
unique separated local proximity space $(X,\rho,\BBBB)$ on
$(X,\tau)$ such that $\BBBB\ap=RC(X)\cap\BBBB$ and
$\rho|_{RC(X)}=\rho\ap$. Moreover, $\BBBB=\{M\sbe X\st \ex
B\in\BBBB\ap$ such that $M\sbe B\}$ and for every $M,N\sbe X$,
$N(-\rho)M\iff M(-\rho)N\iff \fa B\in\BBBB\ \ex F\in\BBBB\ap$ and
$\ex G\in RC(X)$ such that $M\cap B\sbe\int_X(F),\ N\sbe\int_X(G)$
and $F(-\rho\ap)G$.
\end{lm}

\doc The proof that $(X,\rho,\BBBB)$ is a separated local
proximity space on $(X,\tau)$ is straightforward. The uniqueness
follows from Lemma \ref{lead1}.  \sqs

\begin{theorem}\label{lead5}
Let $(X,\tau)$ be a Tychonoff space. Then the ordered sets
$(\LL(X,\tau),\le)$ and $(\LL_{ad}(X,\tau),\preceq_l)$ are
isomorphic.
\end{theorem}

\doc It follows from Lemmas \ref{lead3} and \ref{lead4} and the
Leader Local Compactification Theorem (\cite{LE2}). \sqs

\begin{rem}\label{lead8}
\rm The preceding theorem strengthen Leader Local Compactification
Theorem (\cite{LE2}) because the relation $\rho$ is now given only
on the subset $RC(X)$ of the power set of $X$ and the same is true
for the boundedness $\BBBB$. Using Remark \ref{remd}, we see also
that it is a generalization of the first part of Theorem
\ref{locset} (i.e. that one concerning the ordered set
$(\LL(X),\le)$). In the Leader's paper \cite{LE2} there is no
analogue of the second part of the Theorem \ref{locset} (i.e. that
one concerning the ordered set $(\LL(X),\le_s)$). As far as we
know, another description of the ordered set $(\LL(X),\le_s)$
(even for an arbitrary Tychonoff space $(X,\tau)$) is given only
in \cite[Theorem 5.9]{Di} where this is done by means of the
introduced there special kind of proximities, called
LC-proximities.
\end{rem}

\begin{defi}\label{lead6}
\rm Let $(X,\tau)$ be a Tychonoff space. An NCA $(RC(X,\tau),C)$
is said to be {\em admissible for} $(X,\tau)$ if it satisfies the
following conditions:

\noindent(AK1) if $F,G\in RC(X)$ and $F\cap G\nes$ then $F C G$
(i.e. $\rho_X\sbe C$);

\noindent(AK2) if $F\in RC(X)$ and $x\in\int_X(F)$ then there
exists $G\in RC(X)$ such that $x\in\int_X(G)$ and $G\ll_C F$.

\noindent The set of all NCAs which are admissible for $(X,\tau)$
will be denoted by $\KK_{ad}(X,\tau)$ (or simply by
$\KK_{ad}(X)$). If $(RC(X),C_i)\in\KK_{ad}(X)$, where $i=1,2$,
then we set $(RC(X),C_1)\preceq_{ad} (RC(X),C_2)$ iff $C_2\sbe
C_1$.
\end{defi}

\begin{cor}\label{lead7}{\rm (\cite{dV2})}
If $X$ is a Tychonoff space then the ordered sets $(\KK(X),\le)$
and $(\KK_{ad}(X),\preceq_{ad})$ are isomorphic.
\end{cor}

\doc It follows immediately from Theorem \ref{lead5}. \sqs

Recall that O. Frink \cite{Fr} introduced the notion of a
Wallman-type \cp and asked whether every Hausdorff \cp of a
Tychonoff space is a Wallman-type compactification. This question
was answered in negative by V. M. Ul'janov \cite{Ul}. We will give
a necessary and sufficient condition for a \cp of a discrete space
to be of a Wallman type (recall that, according to the
   Reduction Theorem of L. B. \v{S}apiro \cite{Sha}
   (see also \cite{StSt}), it is
enough to investigate the Frink's problem only in the class of
discrete spaces). Our criterion follows easily from
 the following result of O. Nj\aa stad:

\begin{theorem}\label{Nja}{\rm (\cite{Nj})} Let $(X,\tau)$ be a
Tychonoff space and  $(cX,c)$ be a \cp of $X$. Let $\d_c$ be the
Efremovi\v{c} proximity on the space $(X,\tau)$ defined as
follows: for any $M,N\sbe X$, $M\d_c N$ iff
$\cl_{cX}(c(M))\cap\cl_{cX}(c(N))\nes$. Then $(cX,c)$ is a
Wallman-type \cp of $(X,\tau)$ if and only if there exists a
family $\BB$ of closed subsets of $(X,\tau)$ which is closed under
finite intersections and satisfies the following two conditions:

\noindent(B1) If $F,G\in\BB$ and $F\cap G= \ems$ then $F(-\d_c)
G$;

\noindent(B2) If $A,B\sbe X$ and $A(-\d_c)B$ then there exist
$F,G\in\BB$ such that $A\sbe F$, $B\sbe G$ and $F\cap G=\ems$.
\end{theorem}

\begin{pro}\label{Wallmantype}
Let $X$ be a discrete space, $\Psi^t(X)=(A,\rho,\BBBB)$ and
$C\in\KK_a(X)$. Then $\mu_c\inv(C)$ (see Corollary \ref{compset}
for $\mu_c$) is a Wallman-type \cp of $X$ iff there exists a
sub-meet-semilattice
 $B$ of $A$ such that:

 \noindent(1) for any $a,b\in B$, $a\rho b$ iff $aCb$,  and

\noindent(2) for any $a,c\in A$, $a\ll_C c$ implies that there
exists $b\in B$ such that $a\le b\ll_C c$.
\end{pro}

\doc Note that condition (2) may be substituted for the following one:
\smallskip

\noindent(2\,$\ap$) for any $a,c\in A$, $a\ll_C c$ implies that
there exist $b_1,b_2\in B$ such that $a\le b_1\le b_2^*\le c$.
\smallskip

 Now, applying Remark \ref{remd} and Theorem \ref{Nja}, we complete the proof. \sqs

\baselineskip = 0.75\normalbaselineskip


\begin{thebibliography}{99}
{\small

\bibitem{CNG}
W. Comfort, S.  Negrepontis,
\newblock  Chain Conditions in Topology,
\newblock Cambridge Univ. Press, Cambridge, 1982.


\bibitem{dV2}
H. de Vries,
\newblock  Compact Spaces and Compactifications, an Algebraic Approach,
\newblock Van Gorcum, The Netherlands, 1962.

\bibitem{Di}
G. Dimov,
\newblock   Regular and other kinds of extensions of topological
spaces,
\newblock  Serdica Math. J.  24 (1998) 99--126.


\bibitem{Di2}
 G. Dimov,
\newblock Some generalizations of Fedorchuk Duality Theorem - I,
\newblock arXiv:0709.4495.


\bibitem{DV2}
G. Dimov, D. Vakarelov,
\newblock  Contact algebras and region-based
theory of space: a proximity approach -- I,
\newblock Fundamenta
Informaticae 74(2-3) (2006) 209--249.



\bibitem{E2}
 R. Engelking,
\newblock  General Topology,
\newblock PWN, Warszawa, 1977.






\bibitem{F2}
 V. V. Fedorchuk,
\newblock  Boolean $\d$-algebras and quasi-open mappings,
\newblock  Sibirsk. Mat. \v{Z}. 14 (5) (1973) 1088--1099 = Siberian Math. J.
14 (1973), 759--767 (1974).


\bibitem{Fr}
O. Frink,
\newblock Compactifications and semi-normal spaces,
\newblock Amer. J. Math. 86 (1964) 602--607.



\bibitem{J2}
P. T. Johnstone,
\newblock Stone Spaces,
\newblock Cambridge Univ. Press, Cambridge, 1982.



\bibitem{LE2}
  S. Leader,
\newblock Local proximity spaces,
\newblock  Math. Annalen 169 (1967) 275--281.



\bibitem{MR2}
J. Mioduszewski, L. Rudolf,
\newblock H-closed and extremally disconected Hausdorff spaces,
\newblock  Dissert. Math. (Rozpr. Mat.) 66 (1969) 1--52.

\bibitem{Nj}
O.  Nj\aa stad,
\newblock On Wallman-type compactifications,
\newblock Math. Zeitschr. 91 (1966) 267--276.


\bibitem{Sha}
L. B.  \v{S}apiro,
\newblock A reduction of the basic problem on  compactifications of Wallman type,
\newblock Dokl. Akad. Nauk SSSR  217 (1974) 38--41 = Soviet Math. Dokl. 15 (1974) 1020--1023.


\bibitem{Si2}
  R. Sikorski,
\newblock Boolean Algebras,
\newblock Springer-Verlag, Berlin, 1964.


\bibitem{Sm2}
J. M. Smirnov, On proximity spaces,  Mat. Sb. 31 (1952) 543--574.

\bibitem{StSt}
A. K. Steiner, E. F.  Steiner,
\newblock On the reduction of the Wallman compactification problem to discrete spaces,
\newblock General Topology and Appls. 7 (1977) 35--37.


\bibitem{V2}
 P. Vop\v{e}nka,
\newblock  On the dimension of compact spaces,
\newblock  Czechosl. Math. J. 8 (1958) 319--327.

\bibitem{Ul}
V. M. Ul'janov,
\newblock Solution of a basic problem on
compactifications of Wallman type,
\newblock Soviet Math. Dokl.
 18 (1977) 567--571.


\bibitem{VDDB2}
 D. Vakarelov, G. Dimov, I. D{\"u}ntsch, B. Bennett,
\newblock  A proximity approach to some region-based theories of
space,
\newblock  J. Applied Non-Classical Logics 12 (3-4) (2002)
527--559.

}
\end{thebibliography}
\end{document}